\documentclass[twoside,12pt,a4paper]{article}
\usepackage{latexsym,amsfonts,amstext,amsbsy,amsmath,amssymb,amscd,graphicx,epsfig,eucal,bbm}
\usepackage{mathrsfs}
\usepackage{pb-diagram}
\usepackage[all]{xy}
\usepackage{makeidx}
\usepackage{indentfirst}
\usepackage{multicol}
\usepackage{graphicx}

\newtheorem{theorem}{Theorem}

\newtheorem{pr}[theorem]{Proposition}

\date{}

\usepackage{fancyhdr}
\begin{document}

\title{\textbf{{GENERALIZED FRACTIONAL HYBRID HAMILTON PONTRYAGIN EQUATIONS}}}
\author{\textbf{Chi\c{s} Oana$^1$, Opri\c{s} Dumitru$^2$}\\
\small $ ^1$ West University of Timi\c{s}oara, Romani\\
\small 4 Vasile P\^{a}rvan Blvd., Timi\c{s}oara, 300223, Romania\\
\small chisoana@yahoo.com\\
\small $ ^2$ West University of Timi\c{s}oara, Romania\\
\small 4 Vasile P\^{a}rvan Blvd., Timi\c{s}oara, 300223, Romania\\
\small opris@math.uvt.ro}

\maketitle

\thispagestyle{fancy}

\textbf{Abstract:} \textit{In this work we present a new approach on studying dynamical systems. Combining the two ways of expressing the uncertainty, using probabilistic theory and credibility theory, we have research the generalized fractional hybrid equations. We have introduced the concepts of generalized fractional Wiener process, generalized fractional Liu process and the combination between those two, generalized fractional hybrid process. Corresponding generalized fractional stochastic, respectively fuzzy, respectively hybrid dynamical systems were defined. We applied the theory for generalized fractional hybrid Hamilton-Pontryagin (HP) equation, generalized fractional Hamiltonian equations. From the general fractional hybrid Hamiltonian equations, fractional Langevin equations were found and numerical simulations were done.}\\


\section{INTRODUCTION}

Fractional theory has applicability in many science fields. This approach presents fractional derivatives, fractional integrals, of any real or complex order. Fractional calculus is used when fractional integration is needed. It is used for studying simple dynamical systems, but it also describes complex physical systems. For example, applications of the fractional calculus can be found in chaotic dynamics, control theory, stochastic modeling, but also
in finance, hydrology, biophysics, physics, astrophysics,
cosmology and so on (\cite{chis}, \cite{el-b1}, \cite{el}, \cite{el1}). But
some other fields have just started to study problems from
fractional point of view. In great fashion is the study of
fractional problems of the calculus of variations and
Euler-Lagrange type equations. There were found Euler-Lagrange
equations with fractional derivatives, and then Klimek found
Euler-Lagrange equations, but with symmetric fractional
derivatives \cite{klimek}. Most famous fractional integral are
Riemann-Liouville, Caputo, Grunwald-Letnikov and most frequently
used is Riemann-Liouville fractional integral. The study of
Euler-Lagrange fractional equations was continued  by Agrawal
\cite{agrawal1} that described these equations using the left,
respectively right fractional derivatives in the Riemann-Liouville
sense. This fractional calculus has some great problems, such as
presence of non-local fractional differential operators, or the
adjoint fractional operator that describes the dynamics is not the
negative of itself, or mathematical calculus may be very hard
because of the complicated Leibniz rule, or the absence of chain
rule, and so on. After O.P. Agrawal's formulation \cite{agrawal1}
of Euler-Lagrange fractional equations, B\u{a}leanu and Avkar
\cite{baleanu} used them in formulating problems with Lagrangians
linear in velocities. Standard multi-variable variational calculus
has also some limitations. But in \cite{udriste} C. Udri\c{s}te
and D. Opri\c{s} showed that these limitations can by broken using
the multi-linear control theory.

Another aspect that we use is the stochastic approach. Stochastic concepts were firstly introduced by J.M. Bismut, in his work from 1981, when stochastic Hamiltonian system was introduced. Since then, there has been a need in finding out tools and algorithms for the study of this kind of systems with uncertainty. Bismut's work was continued by Lazaro-Cami and Ortega (\cite{lazaro}, \cite{lazaro1}), in the sense that his work was generalized to manifolds, stochastic Hamiltonian systems on manifolds extremize a stochastic action on the space of manifold valued semimartingales, the reduction of stochastic Hamiltonian system on cotangent bundle of a Lie group, a counter example for the converse of Bismut's original theorem.

A new way for expressing the uncertainty is given in credibility theory. In this case we are not working on a probability space, like in the stochastic case, but on a credibility one. Credibility theory is based on five axioms from which the notion of credibility measure is defined, and it was introduced in order to measure
a fuzzy event. This was first given by Li and Liu in their work \cite{li}. This is a new theory that deals with fuzzy phenomena. Fuzzy random theory and random fuzzy theory can be seen as an extensions of credibility theory. A fuzzy random variable can be seen as a function from a probability space to the set of fuzzy variables, and a random fuzzy variable is a function from a credibility space to the set of random variables \cite{liu1}.  
In our actual research, we will use fuzzy differential equations, that were firstly proposed by Liu \cite{liu}. This is a type of differential equation, driven by a Liu process, just like a stochastic process is described by a Brownian motion.  

In the case when fuzziness and randomness
simultaneously appear in a system, we will talk about hybrid process.  In this sense, we have the concept of fuzzy random variable was introduced by Kwakernaak (\cite{kw1}, \cite{kw2}). A random variable is a random variable that takes fuzzy variable values.  More generally, hybrid variable was proposed by Liu \cite{liu1} to describe the phenomena with fuzziness and randomness. Based on the hybrid process, we will work with differential equations characterized by 
Wiener-Liu process. This can computed using It\^o-Liu formula \cite{you}. In some situations, there exist many Brownian motion (Wiener process) and Liu process in a system, therefore, we can take into consideration also multi-dimensional It\^o-Liu formula.

In this paper, we restrict our attention to stochastic fractional
Hamiltonian systems characterized by Wiener processes and assume
that the space of admissible curves in configuration space is of
class $\mathrm{C}^1.$ Random effects appear in the balance of
momentum equations, as white noise, that is why we may consider
randomly perturbed mechanical systems. It should be mentioned that
the ideas in this paper can be readily extended to stochastic
Hamiltonian systems \cite{milstein} driven by more general
semimartingales, but for the sake of clarity we restrict to Wiener
processes. 

The paper is organized as follows. In Section 2 we present generalization of Riemann-Liouville fractional integral, Wiener process,  and we have defined the generalized fractional stochastic equations.
To get to a hybrid process, we have defined a generalized Liu process and the generalization of fractional fuzzy equation. The mixture between generalized Wiener process and generalized Liu process results as the generalization of hybrid fractional differential equations. In the third section we used the notion presented in Section 2 for defining the generalized fractional hybrid HP equations. We have defined the generalized fractional Riemann-Liouville, respective It\^o, respective Liu integrals and in Theorem 1 we gave generalized fractional hybrid HP equations. We have also defined generalized fractional hybrid Hamiltonian equations. The first order Euler scheme is presented and implemented for particular parameters.

\section{GENERALIZED FRACTIONAL HYBRID EQUATIONS}

Let $f:\mathbb{R}\rightarrow \mathbb{R}$ be an integrable function, $\alpha:\mathbb{R}\rightarrow \mathbb{R}$ a  $C^1$ function, $\rho\in \mathbb{R}, \, \rho>0.$

A \textit{Riemann-Liouville generalized fractional integral} \cite{el1} is defined by
\begin{equation}\label{G1}
_{t_0}I_t^\alpha f(t)=\int_{t_0}^t \frac{1}{\Gamma(\alpha(s-t))}f(s)(t-s)^{\alpha(s-t)-1}e^{-\rho(s-t)}ds,
\end{equation} and 
\begin{equation}\label{G2}
\Gamma_1(\alpha(s-t))=\Gamma(\alpha(z))|_{z=s-t},
\end{equation}and $\Gamma(\alpha(z))$ is Euler Gamma function given by 
\begin{equation}\label{G3} 
\Gamma(\alpha(z))=\int_0^\infty (s-t)^{\alpha(z)-1}e^{-(s-t)}dt.
\end{equation}

If $\alpha(z)=a=const,$ $0<a\leq 1, \, \rho=0,$ from (\ref{G1}), results that
\begin{equation}\label{G4}
_{t_0}I_t^a f(t)=\frac{1}{\Gamma(a)}\int_{t_0}^t f(s)(t-s)^{a-1}ds.
\end{equation}

Formula (\ref{G4}) is the fractional Riemann-Liouville integral, \cite{federico}.

Generalized fractional Riemann-Liouville is a mixture between a fractal action used in physical theory and discount action with rate $\rho,$ \cite{el1}.

In the relations (\ref{G1}) and (\ref{G4}), $s$ is called \textit{intrinsic time} and $t$ is called \textit{observed time}, $t\neq s.$

From (\ref{G1}), results that
\begin{equation}\label{G5}
_{t_0}I_{t}^\alpha f(t)=\int_{t_0}^t f(s)g_t^\alpha(s)ds,
\end{equation} where
\begin{equation}\label{G6}
g_t^\alpha(s)=\frac{1}{\Gamma_1(\alpha(s-t))}e^{(\alpha(s-t)-1)ln|t-s|-\rho(s-t)}, \quad t\neq 0.
\end{equation}

Let $(\Omega, \mathcal{F},P_r)$ be a probabilistic space characterized by the usual conditions, and $(W(t))_{t\in \mathbb{R}}$ a 1-dimensional Wiener process.

It is called a \textit{generalized fractional Wiener process}, the process
\begin{equation}\label{G7}
J^\alpha(t)=\int_{t_0}^t g_t^\alpha(s)dW(s), \quad t\neq s,
\end{equation} where $g_t^\alpha$ is the function given in (\ref{G6}).

If $x(t)=x(t,\omega)$ is a stochastic n-dimensional process and $a:\mathbb{R}\times \mathbb{R}^n\rightarrow \mathbb{R}^n,$
$b:\mathbb{R}\times \mathbb{R}^n\rightarrow \mathbb{R}^n,$ are deterministic functions, we will denote by 
\begin{equation}\label{G8}
\begin{array}{ll}
_{t_0}I_t^\alpha a(t,x(t))=\int_{t_0}^t a(s,x(s))g_t^\alpha (s)ds,\\
   \\
_{t_0}J_t^\beta b(t,x(t))=\int_{t_0}^t b(s,x(s))g_t^\beta (s)dW(s),\\
\end{array}
\end{equation}where $\alpha: \mathbb{R}\rightarrow  \mathbb{R},$
$\beta: \mathbb{R}\rightarrow \mathbb{R},$ $C^1$ functions, the generalized fractional Riemann-Liouville integral, respectively generalized fractional It\^o integral. 

We call \textit{generalized fractional stochastic differential equation}, the functional Volterra type equation given by
\begin{equation}\label{G9}
x(t)=x(t_0)+\int_{t_0}^t a(s,x(s))g_t^\alpha (s)ds+\int_{t_0}^tb(s,x(s))g_t^\beta(s)dW(s).
\end{equation}

Using the notations given (\ref{G8}), it results that 
\begin{equation}\label{G10}
x(t)=x(t_0)+ _{t_0}I_t^\alpha a(t,x(t))+ _{t_0}J_t^\beta b(t,x(t)).
\end{equation}
The equation (\ref{G10}) can be written formally in the following way
\begin{equation}\label{G11}
dx=a(s,x(s))g_t^\alpha(s)ds+b(s,x(s))g_t^\beta(s)dW(s).
\end{equation}

Let us consider $a(t,x(t))=\mu(t)x(t), \, b(t,x(t))=\sigma(t)x(t),$ where $\mu:\mathbb{R}\rightarrow  \mathbb{R}, \ \sigma:\mathbb{R}\rightarrow  \mathbb{R}, \, x:\mathbb{R}\rightarrow  \mathbb{R}$ and $\alpha(z)=\alpha_1, \, \beta(z)=\frac{1+\alpha_1}{2}.$ Then equation (\ref{G9}) becomes 
\begin{equation}\label{4.27}
x(t)=x(t_0)+\frac{1}{\Gamma(\alpha_1)}\int_{t_0}^t \frac{\mu(s)x(s)}{(t-s)^{1-\alpha_1}}ds+\frac{1}{\Gamma(\frac{1+\alpha_1}{2})}\int_{t_0}^t\frac{\sigma(s)x(s)}{(t-s)^{(1-\alpha_1)/2}}ds.
\end{equation} 
The equation (\ref{4.27}) is called \textit{fractional differential equations that governs the stock model (Black-Scholes),} (\cite{black}, \cite{el-b1}, \cite{liu}, \cite{pardoux}, \cite{qin}).

Let $(\Theta, \mathcal{P},C_r)$ be the credibility space with the usual conditions and $(L_t)_{t\in \mathbb{R}}$ an 1-dimensional  Liu process \cite{liu}.

We call a \textit{generalized Liu process}, the following process
\begin{equation}\label{G13}
K^\alpha(t)=\int_{t_0}^t g_t^\alpha(s)dL(s), \quad t\neq s.
\end{equation}

If $x(t)=x(t,\theta)$ is an n-dimensional fuzzy process and
$a_1:\mathbb{R}\times \mathbb{R}^n\rightarrow \mathbb{R}^n,$
$b_1:\mathbb{R}\times \mathbb{R}^n\rightarrow \mathbb{R}^n,$ are deterministic functions, we will denote by
\begin{equation}
\begin{array}{ll}
_{t_0}H_t^\alpha a_1(t,x(t))=\int_{t_0}^t a_1(s,x(s))g_t^\alpha (s)d(s),\\
   \\
_{t_0}L_t^\beta b_1(t,x(t))=\int_{t_0}^t b_1(s,x(s))g_t^\beta (s)dL(s),\\
\end{array}
\end{equation}where $\alpha: \mathbb{R}\rightarrow  \mathbb{R},$
$\beta: \mathbb{R}\rightarrow \mathbb{R}$  are $C^1$ functions, and $g_t^\alpha, \, g_t^\beta$ are given by (\ref{G6}),  the generalized Riemann-Liouville integral, respectively the generalized Liu integral. 

We call \textit{generalized fractional fuzzy differential equation}, the functional Volterra type equation given by
\begin{equation}\label{G15}
\begin{array}{ll}
x(t)=x(t_0)+\int_{t_0}^t a_1(s,x(s))g_t^\alpha (s)ds+\int_{t_0}^tb_1(s,x(s))g_t^\beta(s)dL(s)\\
\quad \quad =x(t_0)+ _{t_0}H_t^\alpha a_1(t,x(t))+ _{t_0}L_t^\beta b_1(t,x(t)).
\end{array}
\end{equation}

Equation (\ref{G15}) can be written formally as
\begin{equation}\label{G11}
dx=a_1(s,x(s))g_t^\alpha(s)ds+b_1(s,x(s))g_t^\beta(s)dL(s).
\end{equation}
 
If $a_1(t,x(t))=\mu(t)x(t), \, b_1(t,x(t))=\sigma(t)x(t), \, \mu:\mathbb{R}\rightarrow\mathbb{R}, \, \sigma:\mathbb{R}\rightarrow\mathbb{R}, \, x:\mathbb{R}\rightarrow\mathbb{R}, \, \alpha(z)=1, \, \beta(z)=\beta_1,$ from (\ref{G10}) results that
\begin{equation}\label{4.31}
x(t)=x(t_0)+\int_{t_0}^t\mu(s)x(s)ds+\frac{1}{\Gamma(\beta_1)}\int_{t_0}^t\frac{\sigma(s)x(s)}{(t-s)^{1-\beta_1}}dL(s).
\end{equation}
The equation (\ref{4.31}) is called \textit{fuzzy equation of a stock model} \cite{qin}.

Let $(\Theta, \mathcal{P},C_r)\times (\Omega, \mathcal{F},P_r)$ be the chance space \cite{liu}, with the usual conditions, and $(W_t)_{t\in \mathbb{R}}$ an 1-dimensional  Wiener process and $(L_t)_{t\in \mathbb{R}}$ an 1-dimensional  Liu process. Let $x(t)=x(t,\omega,\theta)$ an n-dimensional hybrid process and  
$a_2:\mathbb{R}\times \mathbb{R}^n\rightarrow \mathbb{R}^n,$
$b_2:\mathbb{R}\times \mathbb{R}^n\rightarrow \mathbb{R}^n,$ 
$c_2:\mathbb{R}\times \mathbb{R}^n\rightarrow \mathbb{R}^n,$
deterministic functions.

It is called a \textit{generalized fractional hybrid differential equation}, the functional Volterra type equation given by
\begin{equation}\label{G18}
\begin{array}{ll}
x(t)=x(t_0)+\int_{t_0}^t a_2(s,x(s))g_t^\alpha (s)ds+\int_{t_0}^tb_2(s,x(s))g_t^\beta(s)dW(s)\\
\quad \quad +\int_{t_0}^tc_2(s,x(s))g_t^\gamma(s)dL(s),
\end{array}
\end{equation}where $\alpha: \mathbb{R}\rightarrow [0,1], \, \beta:\mathbb{R}\rightarrow \mathbb{R}, \, \gamma:\mathbb{R}\rightarrow \mathbb{R},$ are $C^1$ functions.

With the notations given in (\ref{G8}) and (\ref{G13}), equation (\ref{G18}) can be written as
\begin{equation}\label{G19}
x(t)=x(t_0)+  _{t_0}I_t^\alpha a_2(t,x(t))+ _{t_0}J_t^\beta b_2(t,x(t))+ _{t_0}K_t^\gamma c_2(t,x(t)).
\end{equation}

Formally,  the equation (\ref{G18}) can be expressed as
\begin{equation}\label{G20}
dx=a_2(s,x(s))g_t^\alpha(s)ds+b_2(s,x(s))g_t^\beta(s)dW(s)+c_2(s,x(s))g_t^\gamma(s)dL(s).
\end{equation}

\section{GENERALIZED FRACTIONAL HYBRID HP EQUATION}

Let $Q$ be the paracompact configuration manifold and $J^1(\mathbb{R},Q)=\mathbb{R}\times TQ,$ $T^*Q$ the associated bundle of $Q.$  Let $(\Omega, \mathcal{P},P)$ be a probability space and $(W(t), \mathcal{F}_t)_{t\in [a,b]},$ where $[a,b]\subset \mathbb{R}, \, W(t)$ is a real-valued Wiener process and $\mathcal{F}_t$ is the filtration generated by the Wiener process \cite{bou}. The HP principle unifies the Hamiltonian and Lagrangian description of a mechanical system. The classical HP integral action will be perturbed using deterministic function $\gamma:Q\rightarrow \mathbb{R}.$

Let $\mathcal{L}:J^1(\mathbb{R},Q)\rightarrow \mathbb{R}$ be a $C^2$ function, called \textit{Lagrangian} for the mechanical system and  $\gamma_1,\gamma_2:Q\rightarrow \mathbb{R}$ two functions of class $C^1.$
It is called \textit{generalized fractional action} of $\mathcal{L},$ with respect to the process $(W(t))_{t\in \mathbb{R}}$ and $(L(t))_{t\in \mathbb{R}},$ the function $\mathcal{A}^\alpha:\Theta\times \Omega\times C(PQ)\rightarrow \mathbb{R}$ defined by 
\begin{equation}\label{G21}
\begin{array}{ll}
\mathcal{A}^\alpha(t,q,v,p)=\int_{a}^b(\mathcal{L}(s,q(s),v(s))+<p(s),\frac{dq}{ds}-v(s)>)g_t^\alpha ds +\int_{a}^b\gamma_1(q(s))g_t^\alpha dW(s)\\
\quad \quad \quad \quad \quad \quad +\int_{a}^b\gamma_2(q(s))g_t^\alpha dL(s).
\end{array}
\end{equation}

The first integral in (\ref{G21}) is called \textit{generalized Riemann-Liouville fractional integral}, the second one is \textit{generalized fractional It\^o integral} and the third one is \textit{generalized fractional Liu integral}. Moreover,
$$C(PQ)=\{(t,q,v,p)\in C^0([a,b],J^1(\mathbb{R},Q)), \, q\in C^1([a,b],\mathbb{R}^n), \, q(a)=q_a, \, q(b)=q_b\},$$
$[a,b]\subset \mathbb{R}, \, q_a,q_b\in \mathbb{R}^n.$

We make the following notations 
$q(t,\theta,\omega)=q(t), \, v(t,\theta,\omega)=v(t),\, p(t,\theta,\omega)=p(t).$

Let $c=(q,v,p)\in C([a,b],q_a,q_b)$ be curves on $J^1(\mathbb{R},Q)$ between $q_a$ and $q_b,$ and $B=(q,v,p,\delta q,\delta v,\delta p)\in C^0([a,b], J^1(\mathbb{R},Q)\times J^1(\mathbb{R},Q))$ such that $\delta q(a)=\delta q(b)=0,$ and $q, \, \delta q$ are of class $C^1.$

Let $(q,v,p)(\cdot, \epsilon)\in \mathcal{C}(J^1(\mathbb{R},Q))$ be a family of curves on $J^1(\mathbb{R},Q)$ such that they are differentiable with respect to $\epsilon.$ The differential of the action $\mathcal{A}^\alpha$
is defined by
$$d\mathcal{A}^\alpha(\delta q,\delta v,\delta p)=\frac{\partial}{\partial \epsilon}\mathcal{A}^\alpha(\omega,\theta, q(t,\epsilon),v(t,\epsilon),p(t,\epsilon))\Big|_{\epsilon=0},$$
where 
\begin{equation}\label{G22}
\begin{array}{ll}
\delta q(t)=\frac{\partial}{\partial \epsilon}q(t,\epsilon)\Big|_{\epsilon=0}, \, \delta q(a)=\delta q(b), \\
   \\
\delta v(t)=\frac{\partial}{\partial \epsilon}v(t,\epsilon)\Big|_{\epsilon=0}, \, 
\delta p(t)=\frac{\partial}{\partial \epsilon}p(t,\epsilon)\Big|_{\epsilon=0}.
\end{array}
\end{equation}
Using (\ref{G21}), by direct calculus, we get the following theorem.

\begin{theorem}
Let $\mathcal{L}:J^1(\mathbb{R},Q)\rightarrow \mathbb{R}$ be a Lagrangian $C^2$  function with respect to $t,\, q$ and $ v$ and the first order derivatives are Lipschitz functions with respect to $t, \, q, \, v.$ Let $\gamma_1,\gamma_2:Q\rightarrow \mathbb{R}$ be functions of class $C^2,$ and with the first order derivatives Lipschitz functions. Then, the curve $c=(q,v,p)\in C(J^1(\mathbb{R},Q)\times \mathbb{R}^n)$ satisfies the  generalized fractional hybrid HP equations a.s.
\begin{equation}\label{G23}
\begin{array}{ll}
dq^i=v^i ds,\\
  \\
dp_i=(\frac{\partial \mathcal{L}}{\partial q^i}-p_ih(s,t))ds+
\frac{\partial \gamma_1(q)}{\partial q^i}dW(s)+\frac{\partial \gamma_2(q)}{\partial q^i}dL(s),\\
  \\
p_i=\frac{\partial \mathcal{L}}{\partial v^i}, \, i=1,...,n, \, t\neq s,
\end{array}
\end{equation}where 
$$h(s,t)=\frac{d (\alpha(s-t))}{ds}ln |t-s|+\frac{\alpha(s-t)-1}{s-t}+\rho-\frac{1}{\Gamma_1(\alpha(s-t))}\frac{d\Gamma_1(\alpha(s-t))}{ds}.$$\hfill $\Box$
\end{theorem}

From (\ref{G23}) we have:
\begin{description}
	\item[(i)] If $\alpha(z)=1,\, \rho=0,$ then
	\begin{equation}\label{G24}
\begin{array}{ll}
dq^i=v^i ds,\\
   \\
dp_i=\frac{\partial \mathcal{L}}{\partial q^i}ds+
\frac{\partial \gamma_1(q)}{\partial q^i}dW(s)+\frac{\partial \gamma_2(q)}{\partial q^i}dL(s),\\
   \\
p_i=\frac{\partial \mathcal{L}}{\partial v^i}, \, i=1,...,n, \, t\neq s;
\end{array}
\end{equation}
	\item[(ii)] If $\alpha(z)=a=const,\, 0<a\leq 1, \,  
	\rho=0,$ then  
    \begin{equation}\label{G25}
\begin{array}{ll}
dq^i=v^i ds,\\
   \\
dp_i=(\frac{\partial \mathcal{L}}{\partial q^i}-p_i\frac{a-1}{s-t})ds+
\frac{\partial \gamma_1(q)}{\partial q^i}dW(s)+\frac{\partial \gamma_2(q)}{\partial q^i}dL(s),\\
  \\
  p_i=\frac{\partial  \mathcal{L}}{\partial v^i}, \, ,\, i=1,...,n, \, t\neq s.
\end{array}
\end{equation}
\end{description}

For $\gamma_2=0,$ we get the fractional stochastic HP 
equations \cite{chis4}.

If $\mathcal{L}:M\rightarrow \mathbb{R}$ is hyperregular, that means $det \Big(\frac{\partial^2 \mathcal{L}}{\partial v^i \partial v^j} \Big)\neq 0,$ from (\ref{G23}) results the following proposition.

\begin{pr}(Generalized fractional hybrid Hamiltonian equations)\\
The equations \emph{(\ref{G23})} are equivalent with the equations
\begin{equation}\label{G26}
\begin{array}{ll}
dq^i=\frac{\partial H}{\partial p_i}ds, \\
   \\
dp_i=(-\frac{\partial H}{\partial q^i}-p_ih(s,t))ds+
\frac{\partial \gamma_1(q)}{\partial q^i}dW(s)+\frac{\partial \gamma_2(q)}{\partial q^i}dL(s), \\
\end{array}
\end{equation}where
$$
H=p_iv^i-\mathcal{L}(t,q,v),$$
$$h(s,t)=\frac{d \alpha(s-t)}{ds}ln|t-s|+\frac{\alpha(s-t)-1}{s-t}+\rho-\frac{1}{\Gamma_1(\alpha(s-t))}\frac{d\Gamma_1(\alpha(s-t))}{ds}$$
\end{pr}\hfill $\Box$

The equations (\ref{G26}) are called \textit{generalized fractional hybrid Langevin equations} and  can be written to describe the movement equations for relativistic particles with white noise and Liu process.

\begin{pr}
If $\mathcal{L}=\frac{1}{2}g_{ij}v^iv^j,$ where $g_{ij}$ are the components of a metric on a manifold $Q,$ then the equations \emph{(\ref{G23})} take the form
\begin{equation}\label{eq3}
\begin{array}{ll}
dq^i=v^ids, \\
  \\
dv^i=-(\Gamma_{jk}^i v^jv^k-h(s,t)v^i)ds+g^{ij}\frac{\partial \gamma_1(q)}{\partial q^j}dW(s)+g^{ij}\frac{\partial \gamma_2(q)}{\partial q^j}dL(s), \, i,j=1,...,n,
\end{array}
\end{equation}where $\Gamma_{jk}^i$ are Cristoffel coefficients associated to the considered metric and $h(s,t)$ is given above.

The equations \emph{(\ref{G26})} become
\begin{equation}\label{s1}
\begin{array}{ll}
dq^i=g^{ij}p_jds, \\
   \\
dp_i=(\frac{1}{2}\frac{\partial g_{kl}}{\partial q^i}p^kp^l-h(s,t)p_i)ds+\frac{\partial \gamma_1(q)}{\partial q^i}dW(s)
+\frac{\partial \gamma_2(q)}{\partial q^i}dL(s), \, i=1,...,n.
\end{array}
\end{equation}\hfill $\Box$
\end{pr}

\begin{pr}
Let $\mathcal{L}:J^1(\mathbb{R},\mathbb{R})\rightarrow \mathbb{R}$ be given by
$$\mathcal{L}(q,v)=\frac{1}{2}v^2-V(q),$$ and $V,\gamma_1,\gamma_2:\mathbb{R}\rightarrow \mathbb{R}.$ The equations \emph{(\ref{s1})} are given by
\begin{equation}\label{G27}
\begin{array}{ll}
dq=pds, \\
   \\
dp=(-\frac{\partial V}{\partial q}-h(s,t)p)ds+
\frac{\partial \gamma_1(q)}{\partial q}dW(s)+\frac{\partial \gamma_2(q)}{\partial q}dL(s).
\end{array}
\end{equation}
\end{pr}\hfill $\Box$

If $V(q)=\cos(q),\, \gamma_1(q)=\alpha_1\sin(q)$ and $\gamma_2(q)=\frac{1}{2}\alpha_2q^2,$ the first order Euler scheme for the equations (\ref{G27}) is given by
\begin{equation}
\begin{array}{ll}
q(n+1)=q(n)+Kp(n),\\
   \\
p(n+1)=p(n)+K(\sin(q(n))-h(nK,t))+\alpha_1\cos(q(n))G(n)+\alpha_2q(n)L(n,z_2),  
\end{array}
\end{equation}where $n=0,...,N-1, \, K=\frac{T}{N},$ $G(n)$ and $L(n,z)$ are the simulations of Wiener and Liu processes and 
$$h(nK,t)=\dot{\alpha}(nK,t)ln|t-nK|+\frac{\alpha(nK-t)-1}{nK-t}+\rho-\frac{1}{\Gamma_1(\alpha(nK-t))}\dot{\Gamma_1}(\alpha(nK-t)),$$ with $$\dot{\alpha}(s,t)=\frac{d\alpha(s-t)}{ds}, \quad \dot{\Gamma_1}(\alpha(s,t))=\frac{d\Gamma_1(\alpha(s-t))}{ds},$$ and 
$$G(n) = random[normald[0, \sqrt{h}](1),$$
$$L(n,z)= \frac{2}{1+e^{\pi |z|/(h\sigma\sqrt{6}S_2(n))}}, \, S_2(n)=\sum_{k=0}^{n-1}(b(\alpha_2 q(k),z(k))).$$

Using Maple 13, for the values of the parameters, $\alpha=0.6, \, t=0.8, \, \alpha_1=0.1, \, \alpha_2=0.3,\, z_2=15.,$ we get the following orbits. 

\begin{center}\begin{tabular}{ccc}
\epsfxsize=5cm \epsfysize=4cm
 \epsffile{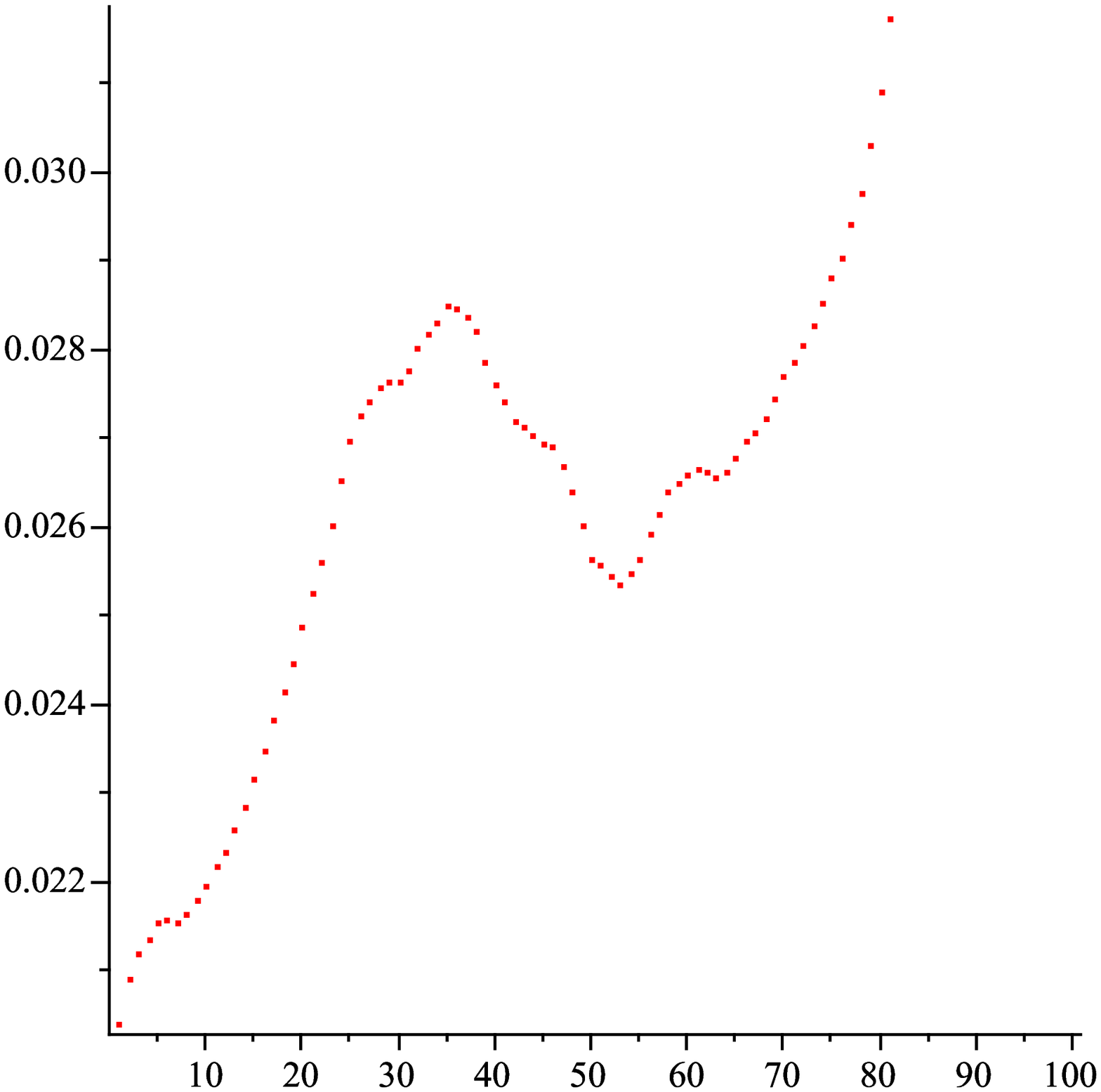}&
 \epsfxsize=5cm \epsfysize=4cm
\epsffile{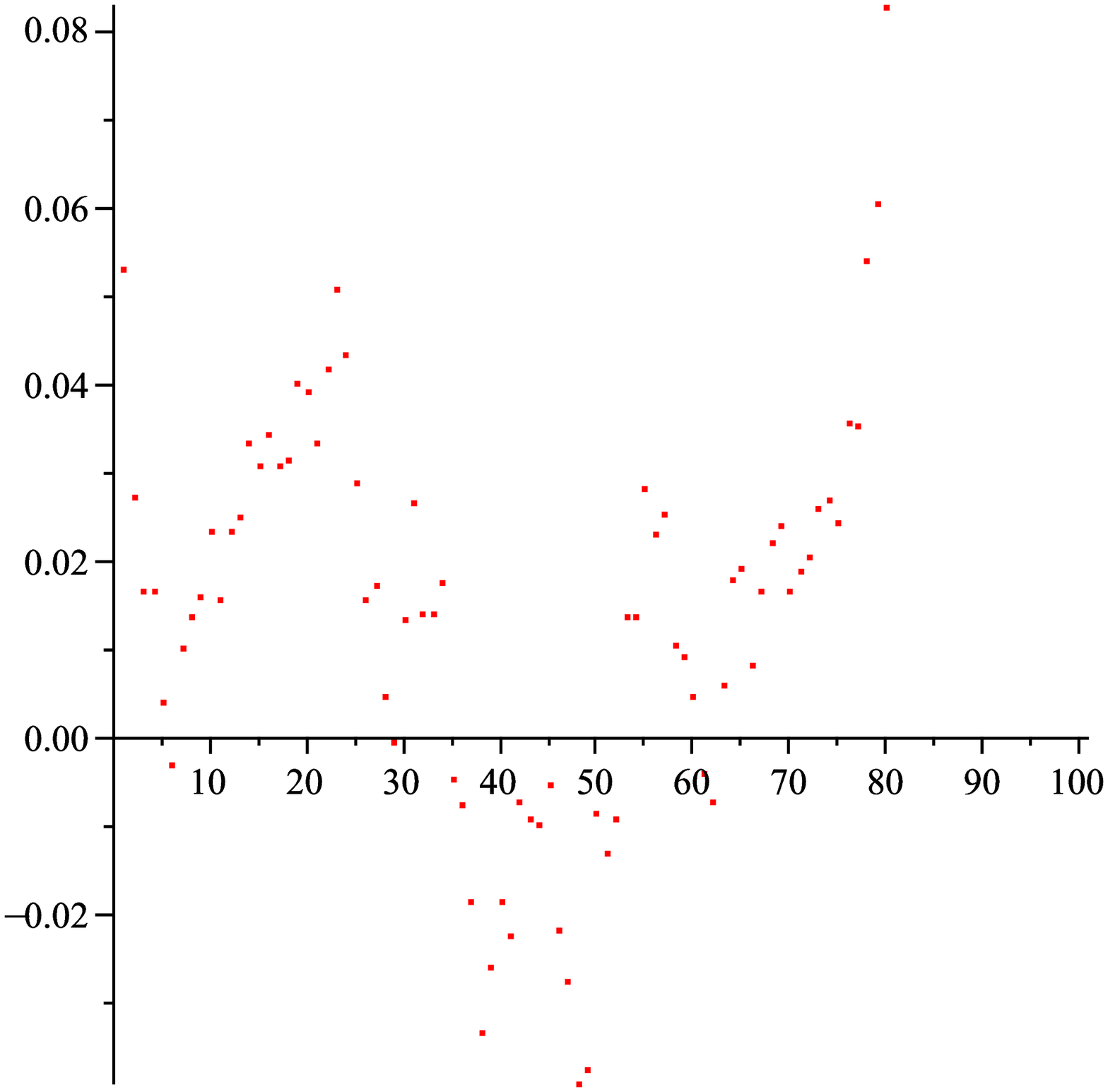}& 
\epsfxsize=5cm \epsfysize=4cm
\epsffile{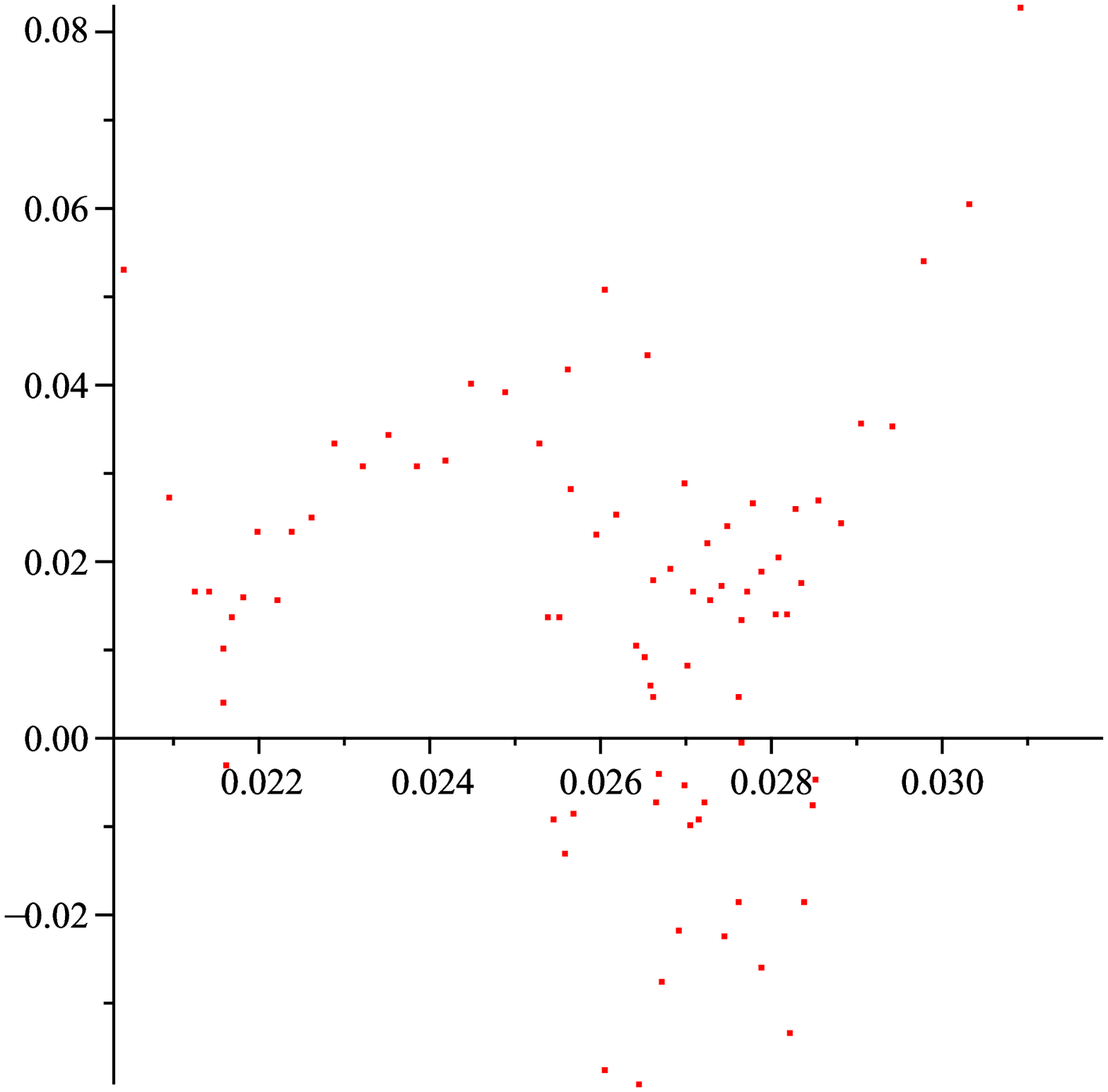}\\
\textit{ Fig1:} $(n,q(n,z_2,\omega))$& \textit{ Fig2: }$(n,p(n,z_2,\omega))$&
\textit{ Fig3: }$(q(n,z_2,\omega),p(n,z_2,\omega))$\\
    \end{tabular}
\end{center}

If we are in the classical case, with $\alpha=1,$ then the above graphics become

\begin{center}\begin{tabular}{ccc}
\epsfxsize=5cm \epsfysize=4cm
 \epsffile{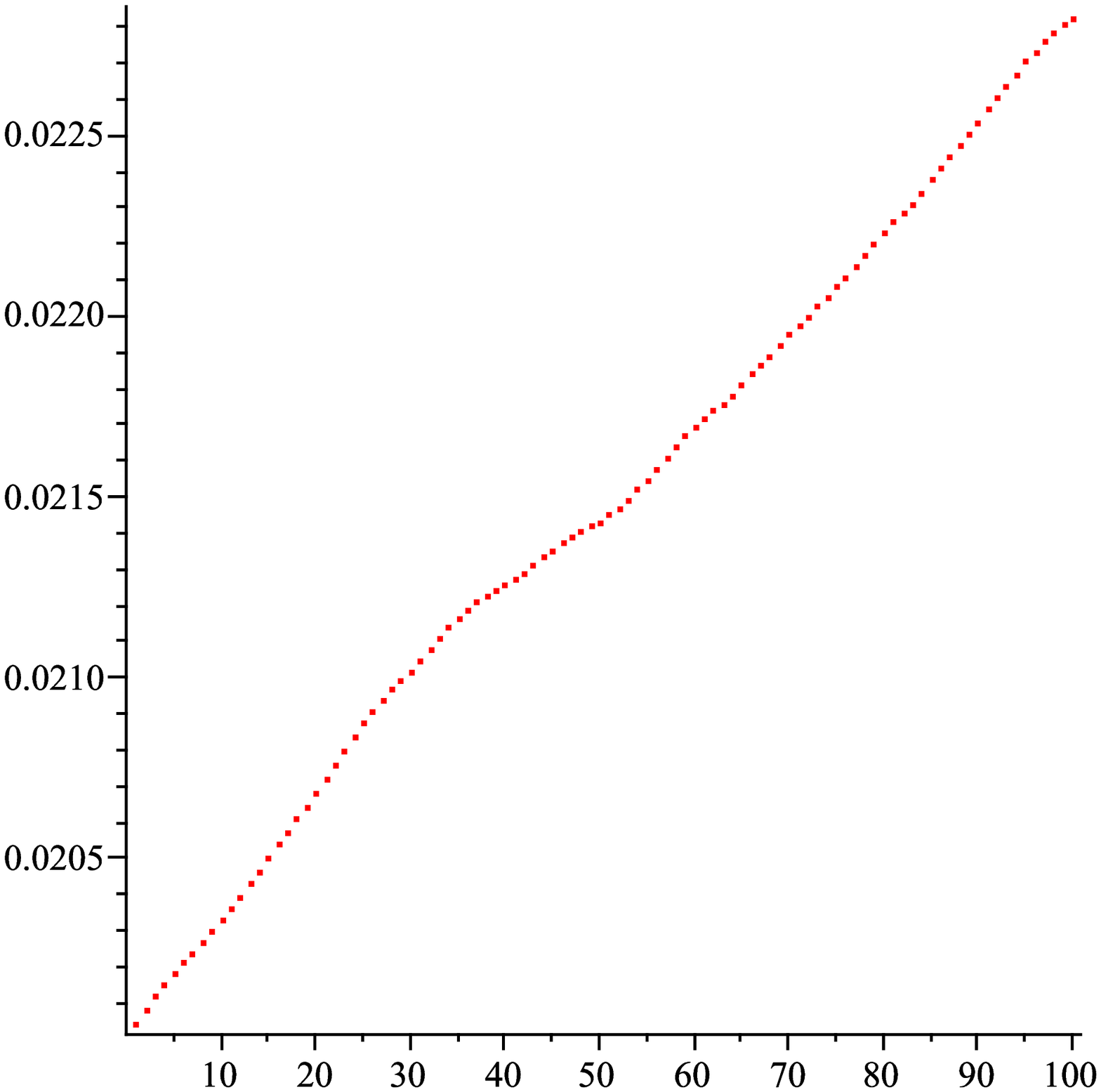}&
 \epsfxsize=5cm \epsfysize=4cm
\epsffile{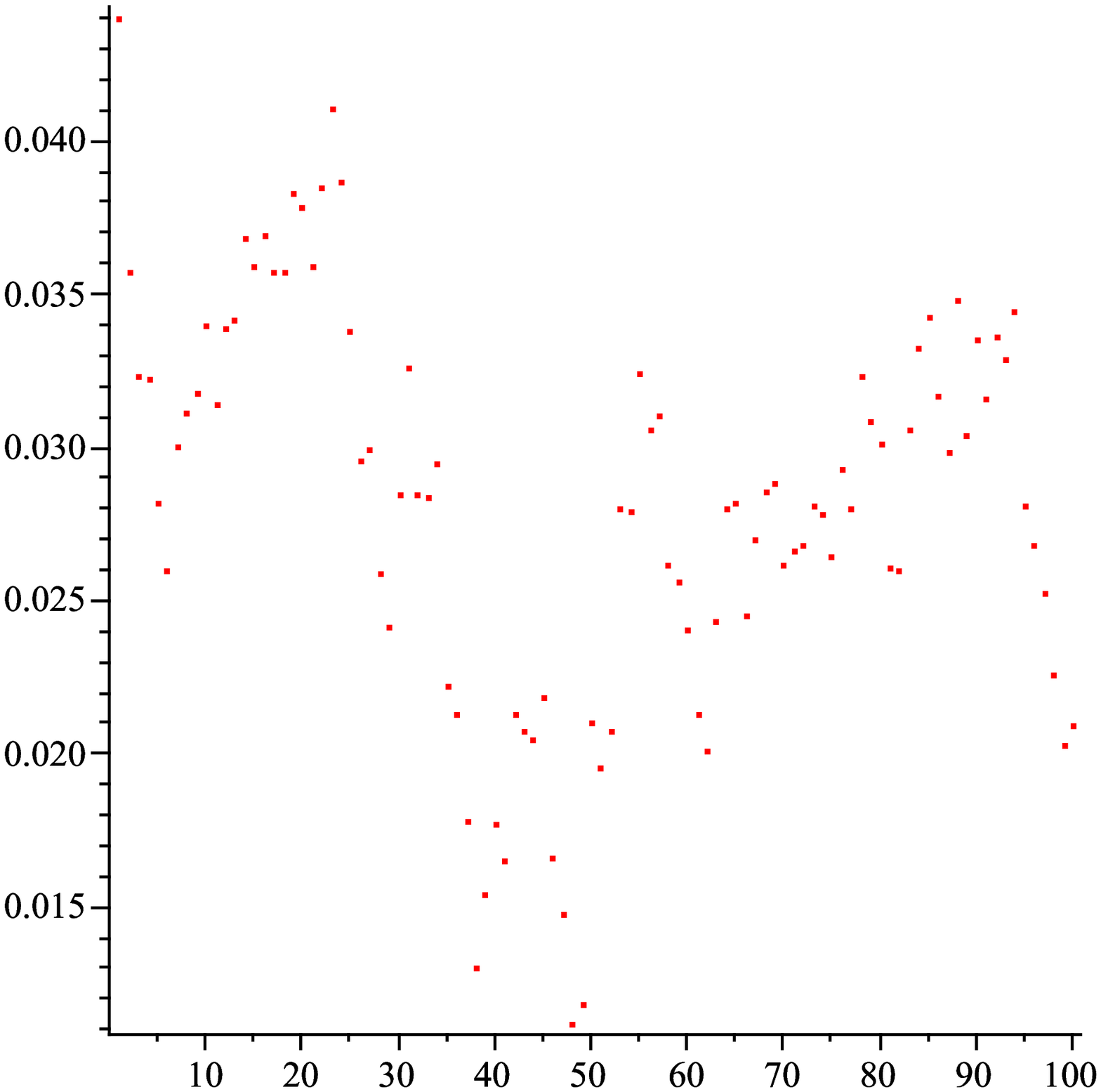}& 
\epsfxsize=5cm \epsfysize=4cm
\epsffile{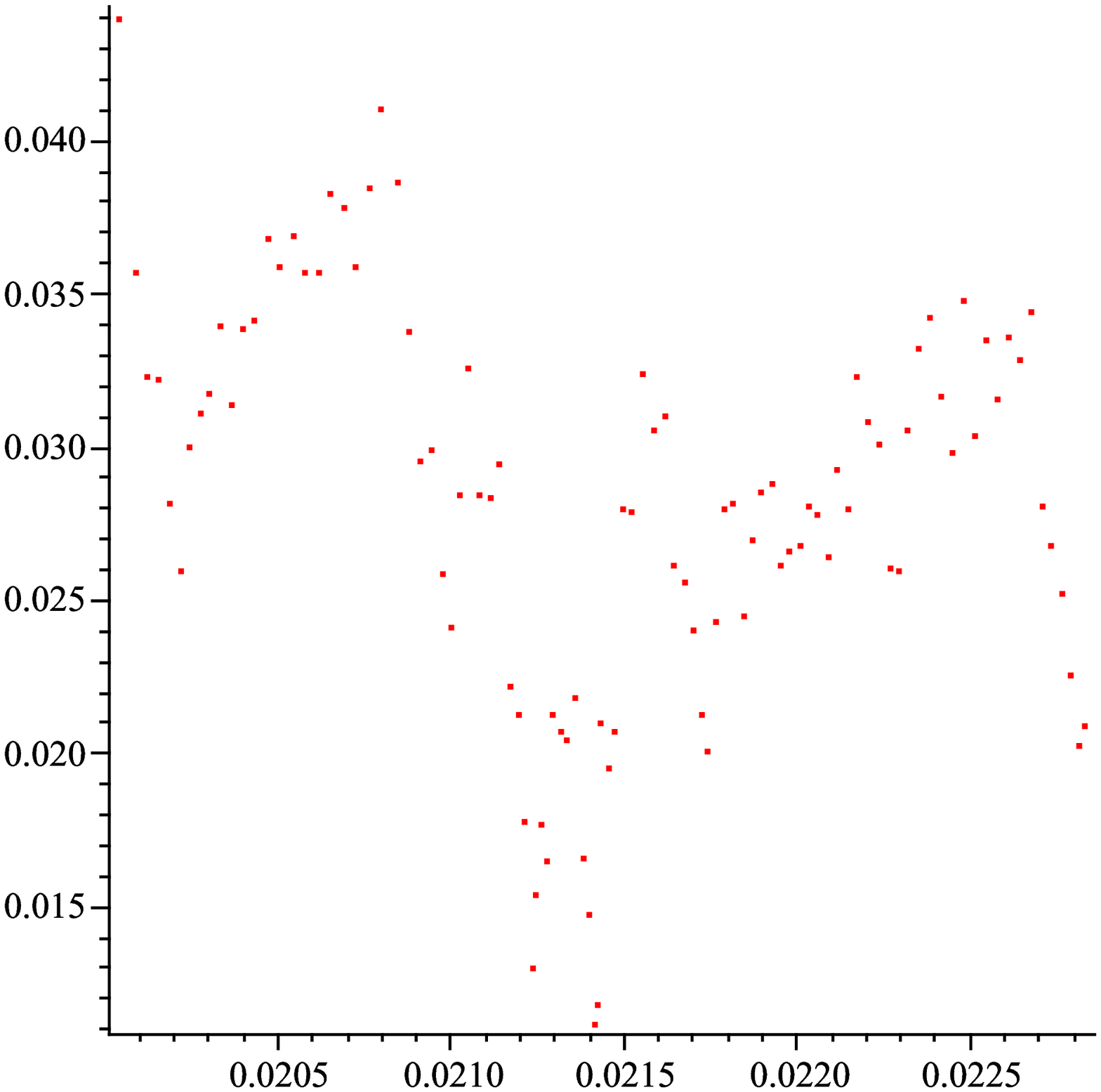}\\
\textit{ Fig4:} $(n,q(n,z_2,\omega))$& \textit{ Fig5:} $(n,p(n,z_2,\omega))$&
\textit{ Fig6:} $(q(n,z_2,\omega),p(n,z_2,\omega))$\\
    \end{tabular}
\end{center}

For $\alpha=0.6, \, t=0.8, \, \alpha_1=0., \, \alpha_2=0.3, \, z_2=15.,$ we get the orbits given in figures 7, 8 and 9, and if $\alpha=1,$ we get the figures 10, 11 and 12. 

\begin{center}\begin{tabular}{ccc}
\epsfxsize=5cm \epsfysize=4cm
 \epsffile{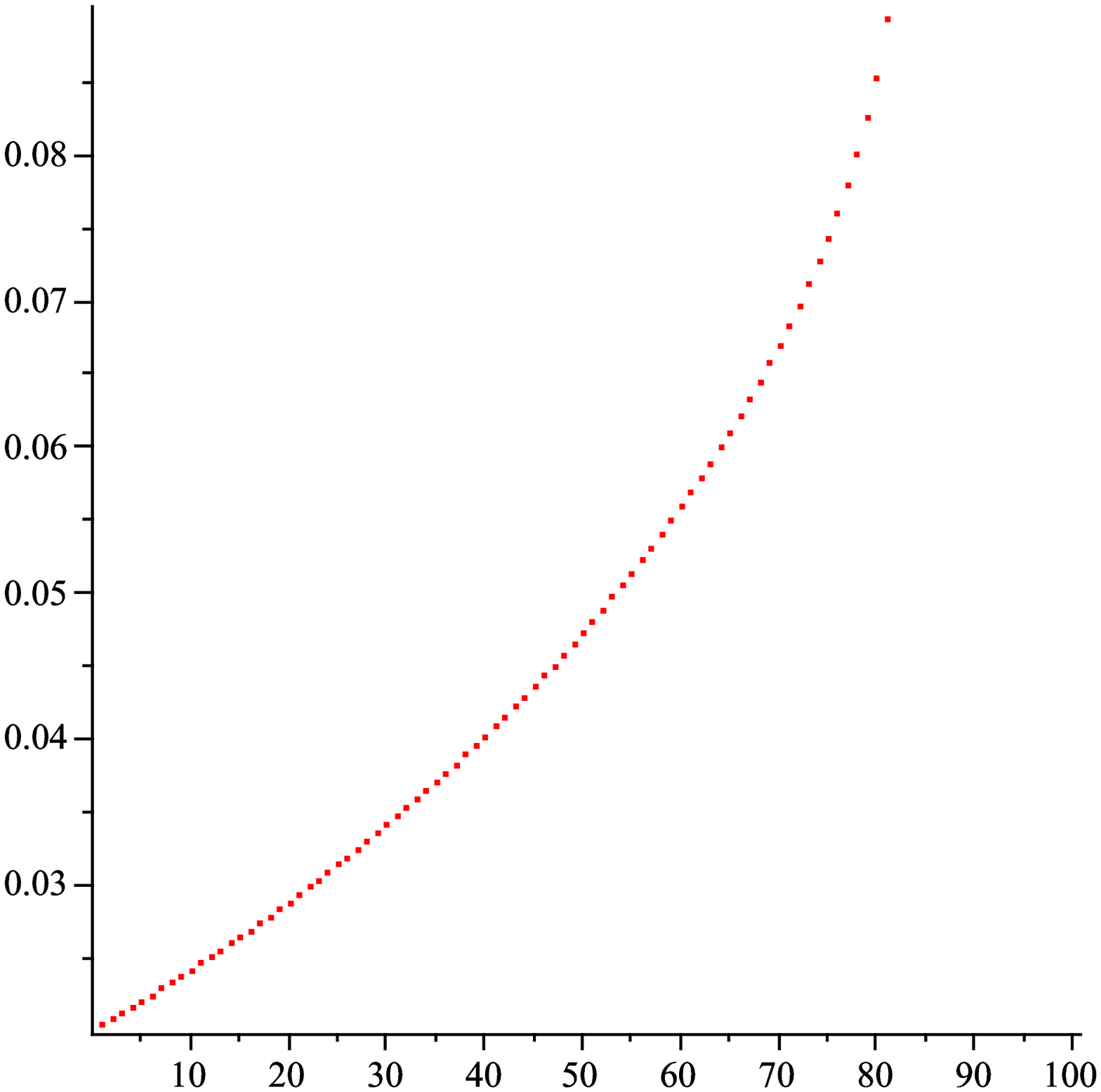}&
 \epsfxsize=5cm \epsfysize=4cm
\epsffile{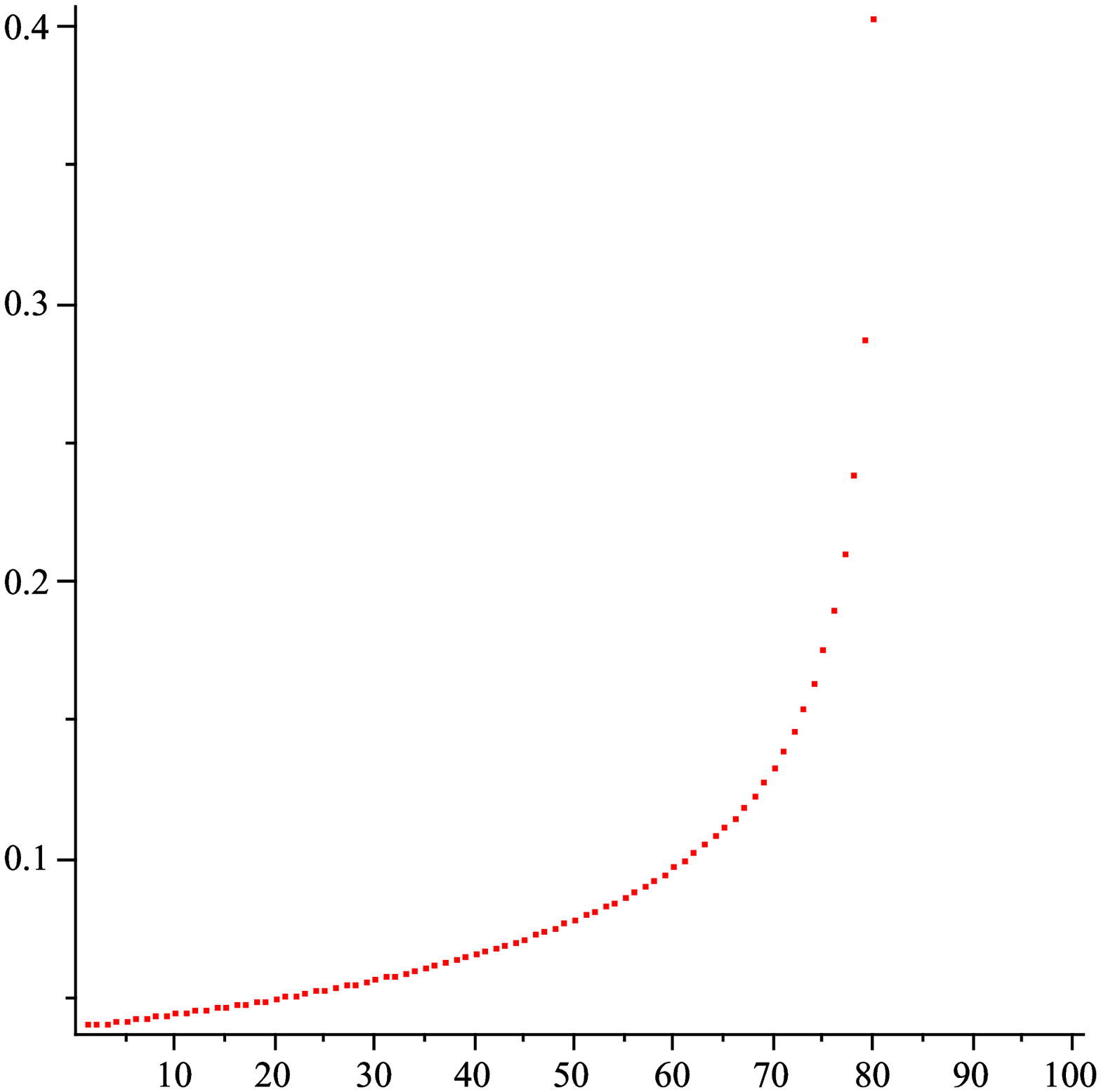}& 
\epsfxsize=5cm \epsfysize=4cm
\epsffile{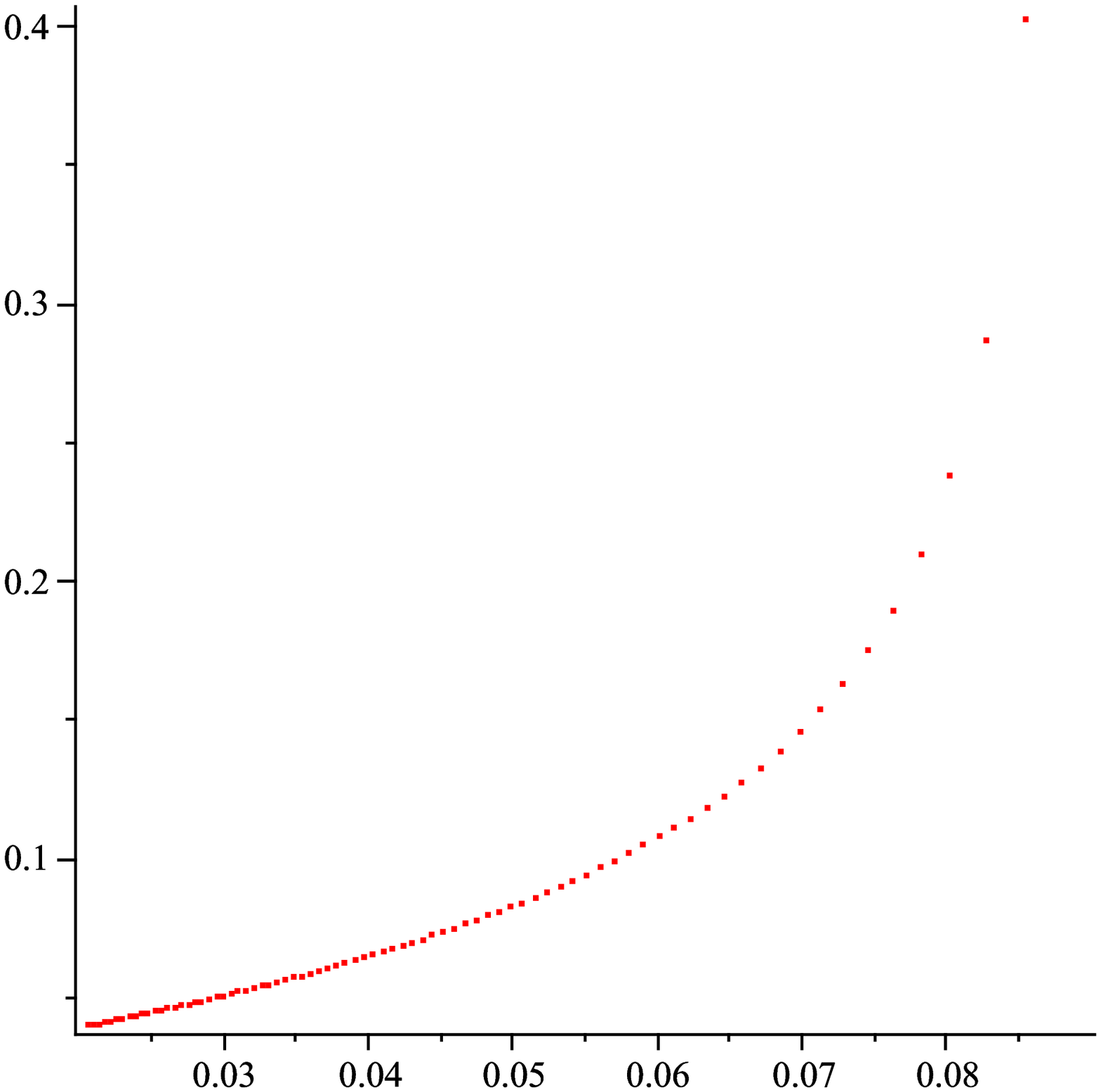}\\
\textit{ Fig7:} $(n,q(n,z_2))$& \textit{ Fig8:} $(n,p(n,z_2))$&
\textit{ Fig9:} $(q(n,z_2),p(n,z_2))$\\
    \end{tabular}
\end{center}

\begin{center}\begin{tabular}{ccc}
\epsfxsize=5cm \epsfysize=4cm
 \epsffile{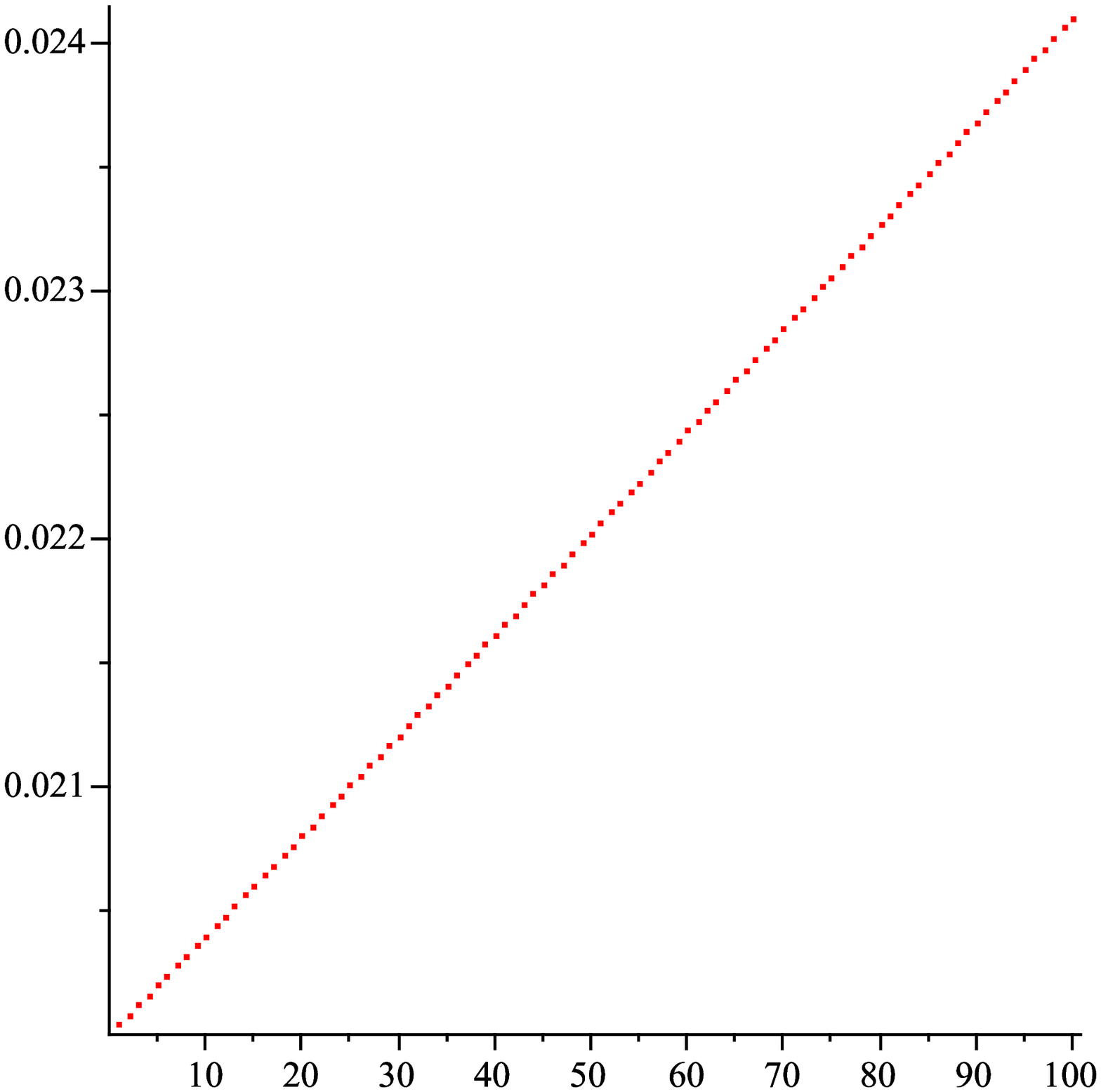}&
 \epsfxsize=5cm \epsfysize=4cm
\epsffile{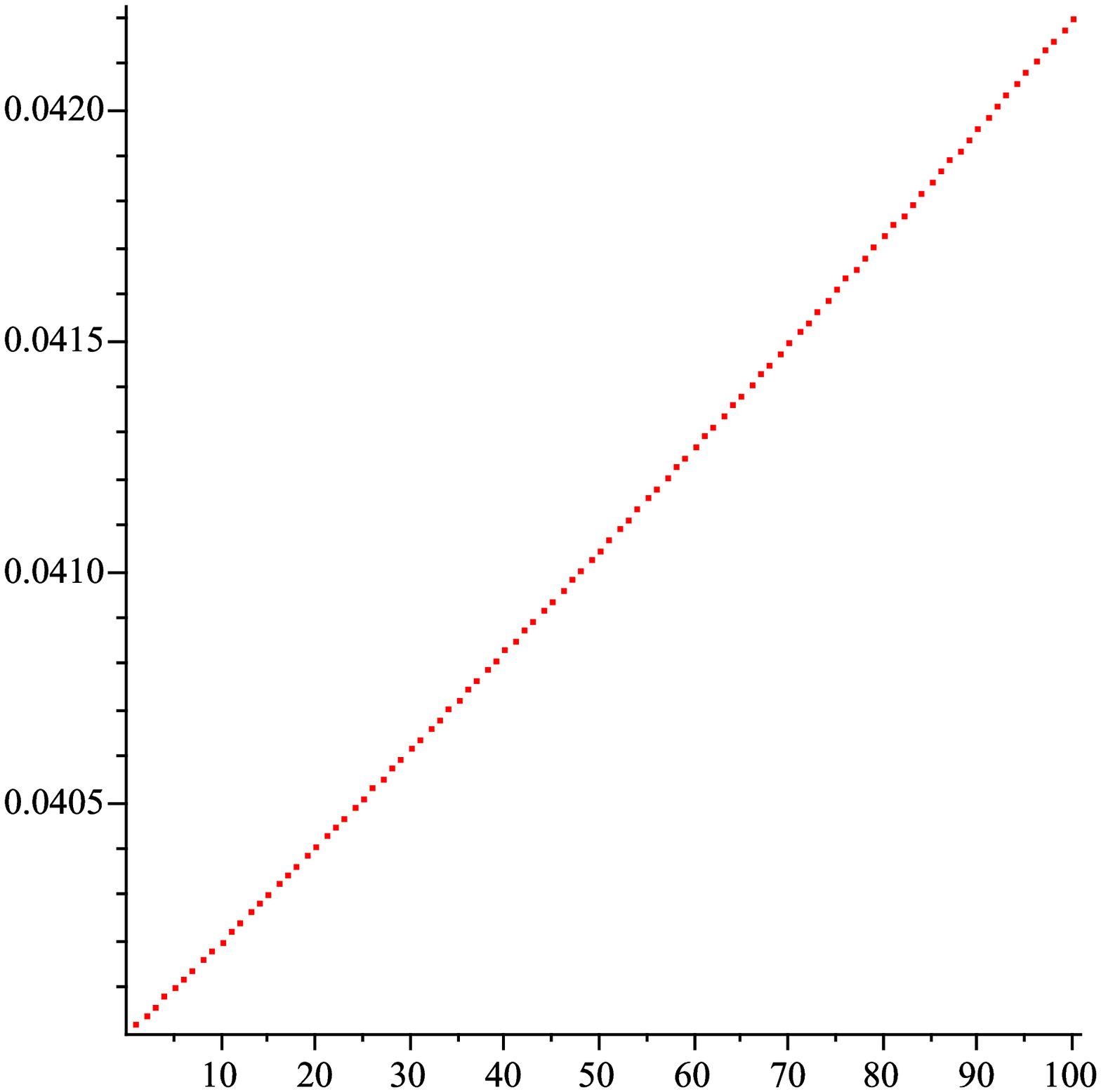}& 
\epsfxsize=5cm \epsfysize=4cm
\epsffile{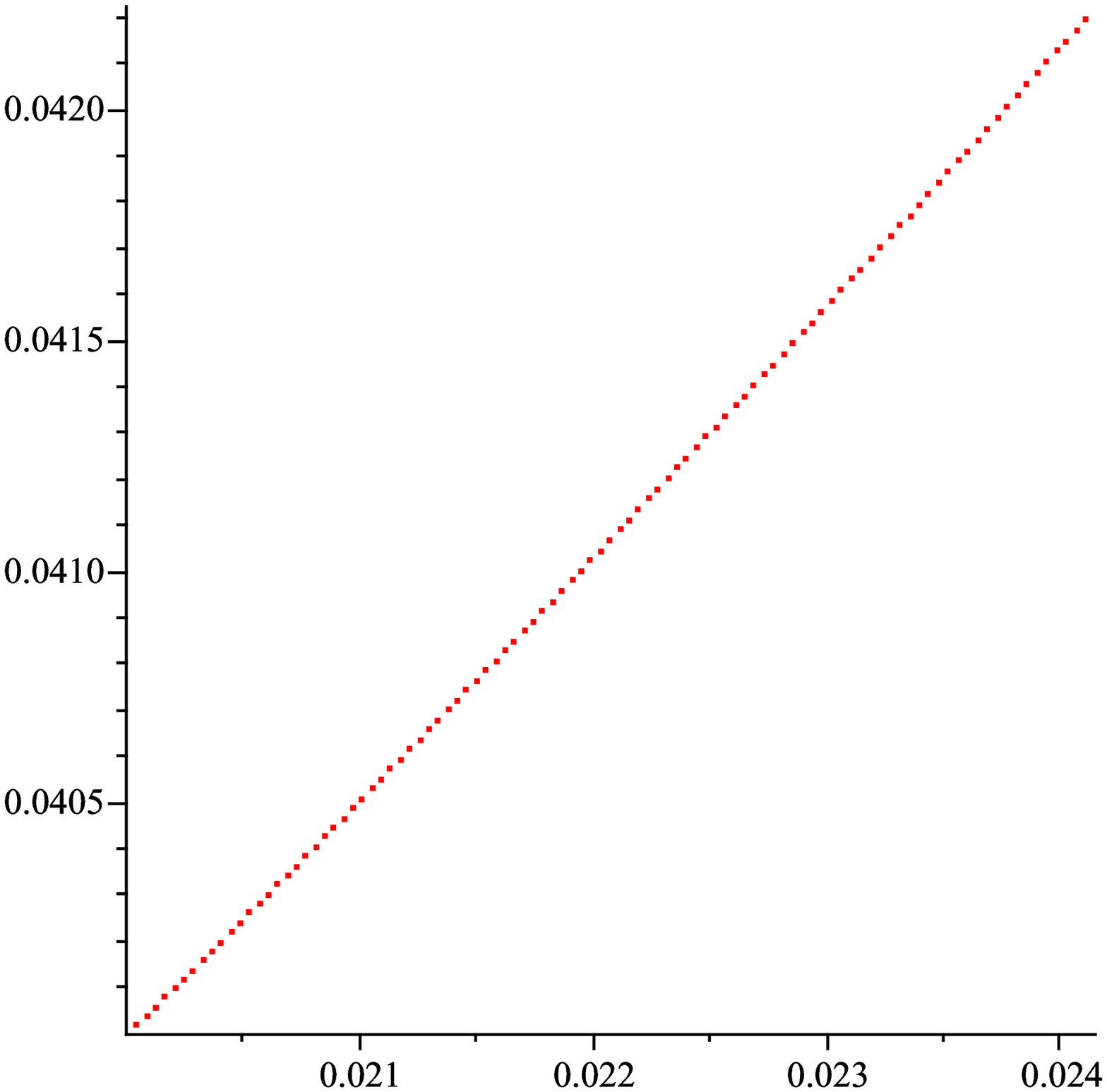}\\
\textit{ Fig10: }$(n,q(n,z_2))$& \textit{ Fig11:} $(n,p(n,z_2))$&
\textit{ Fig12:} $(q(n,z_2),p(n,z_2))$\\
    \end{tabular}
\end{center}

The orbits for $(n,q(n,\omega)), \, (n,p(n,\omega)), \, (q(n,\omega),p(n,\omega)),$ for the values of the parameters $\alpha=0.6, \, t=0.8, \, \alpha_1=0.1, \, \alpha_2=0, \, z_2=15.,$are represented in figures 13, 14 and 15.

\begin{center}\begin{tabular}{ccc}
\epsfxsize=5cm \epsfysize=4cm
 \epsffile{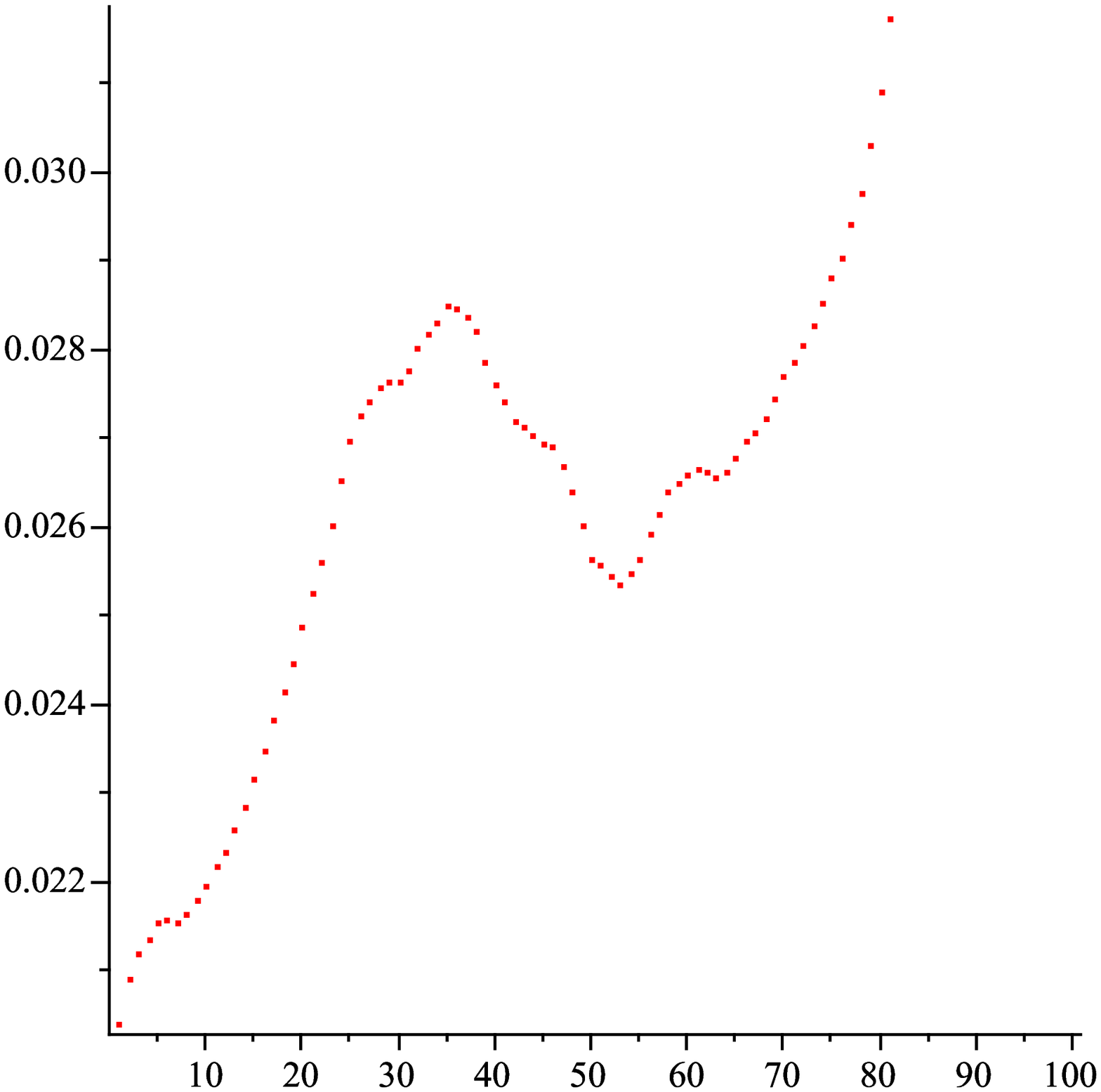}&
 \epsfxsize=5cm \epsfysize=4cm
\epsffile{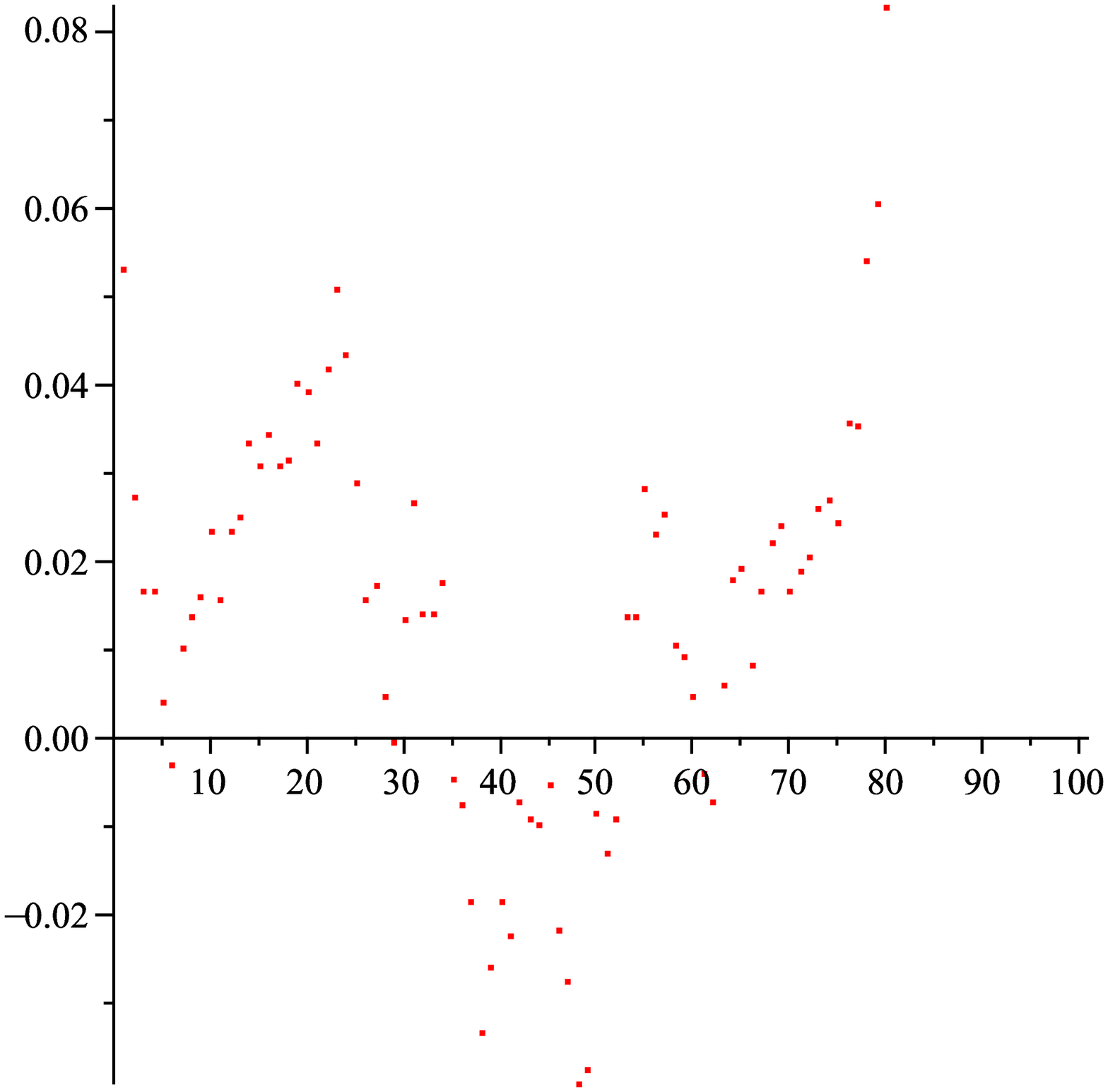}& 
\epsfxsize=5cm \epsfysize=4cm
\epsffile{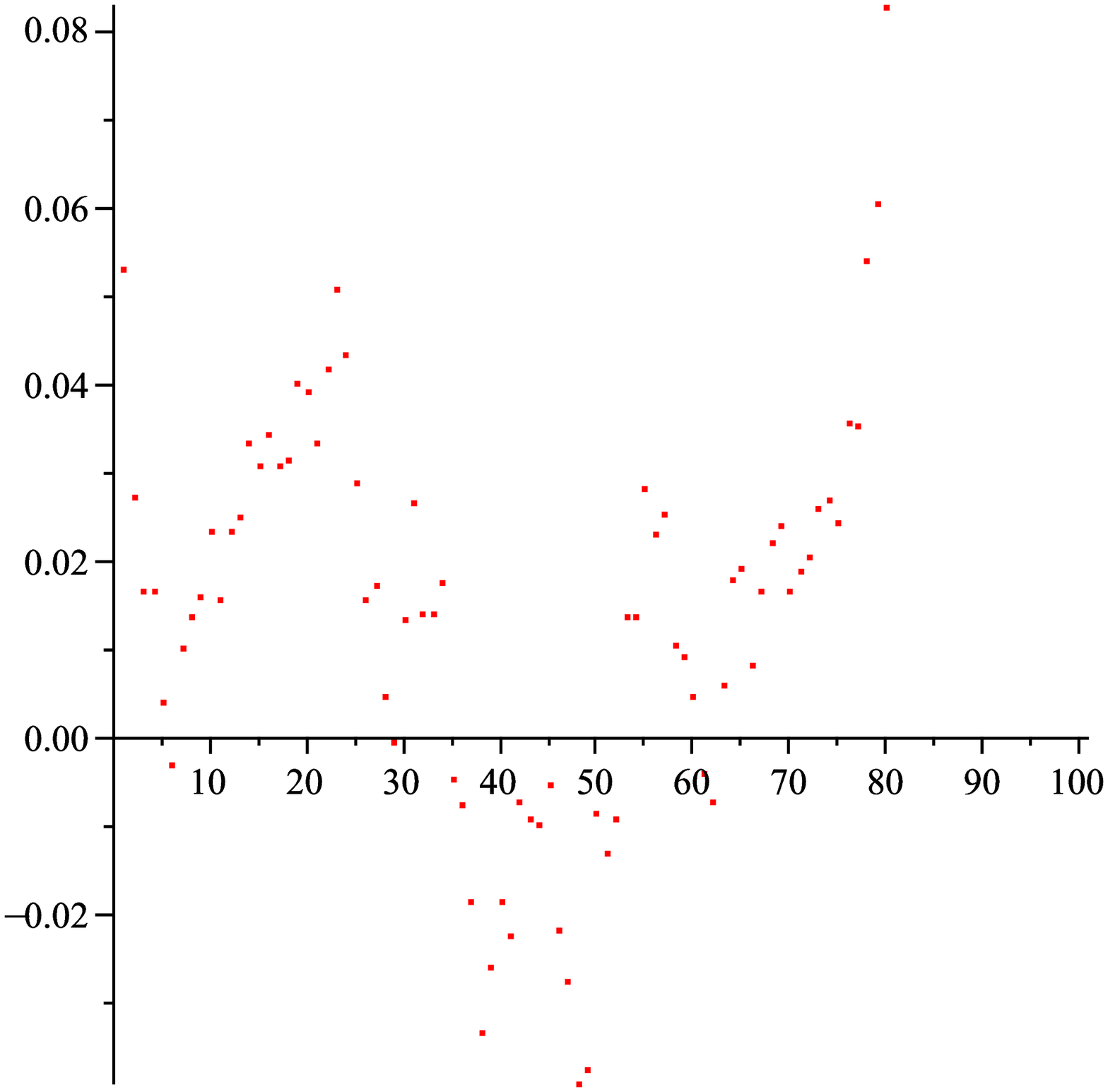}\\
\textit{ Fig13:} $(n,q(n,\omega))$& \textit{ Fig14:} $(n,p(n,\omega))$&
\textit{ Fig15: }$(q(n,\omega),p(n,\omega))$\\
    \end{tabular}
\end{center}

For $\alpha=1,$ the figures 13, 14 and 15 become

\begin{center}\begin{tabular}{ccc}
\epsfxsize=5cm \epsfysize=4cm
 \epsffile{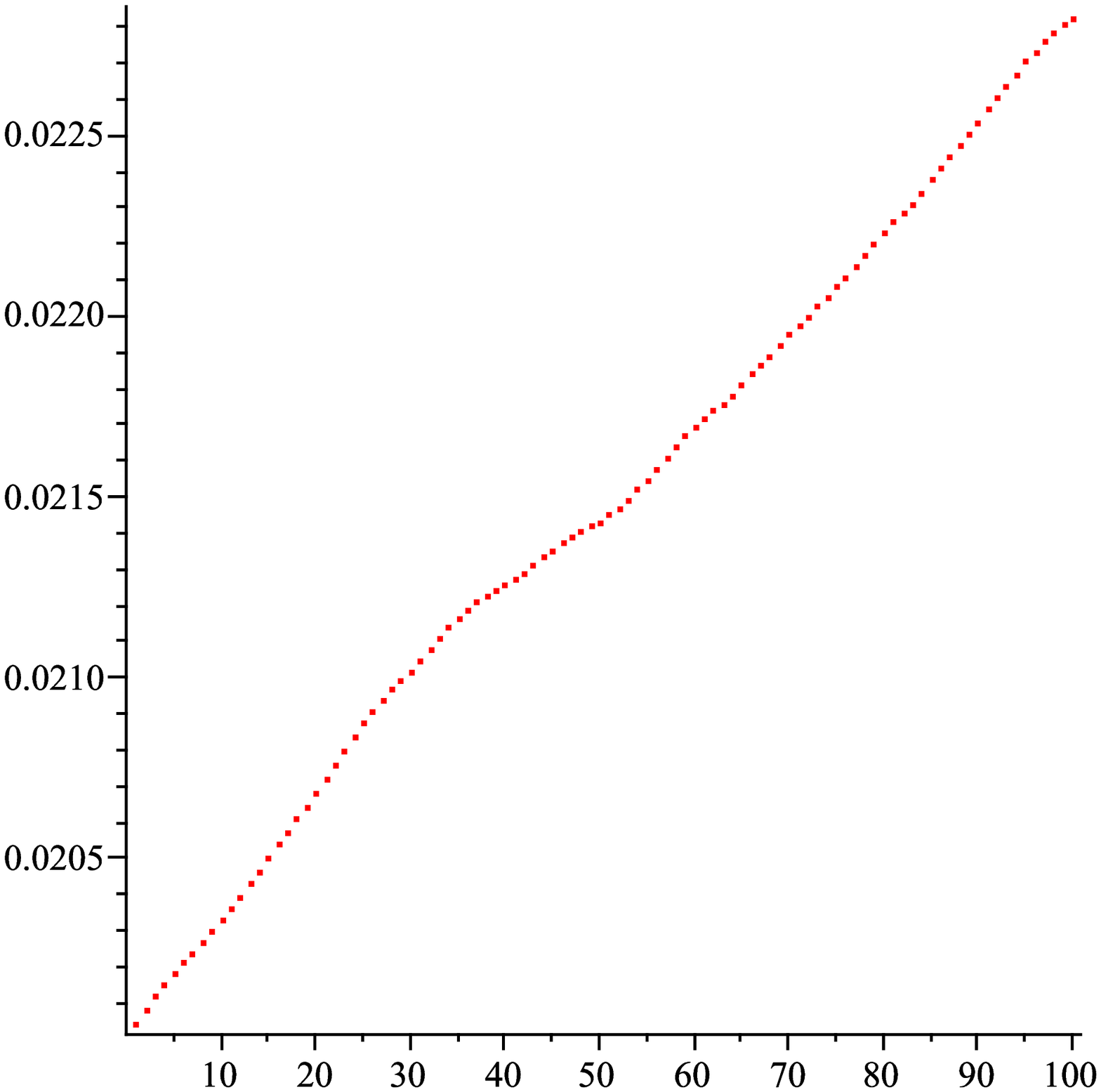}&
 \epsfxsize=5cm \epsfysize=4cm
\epsffile{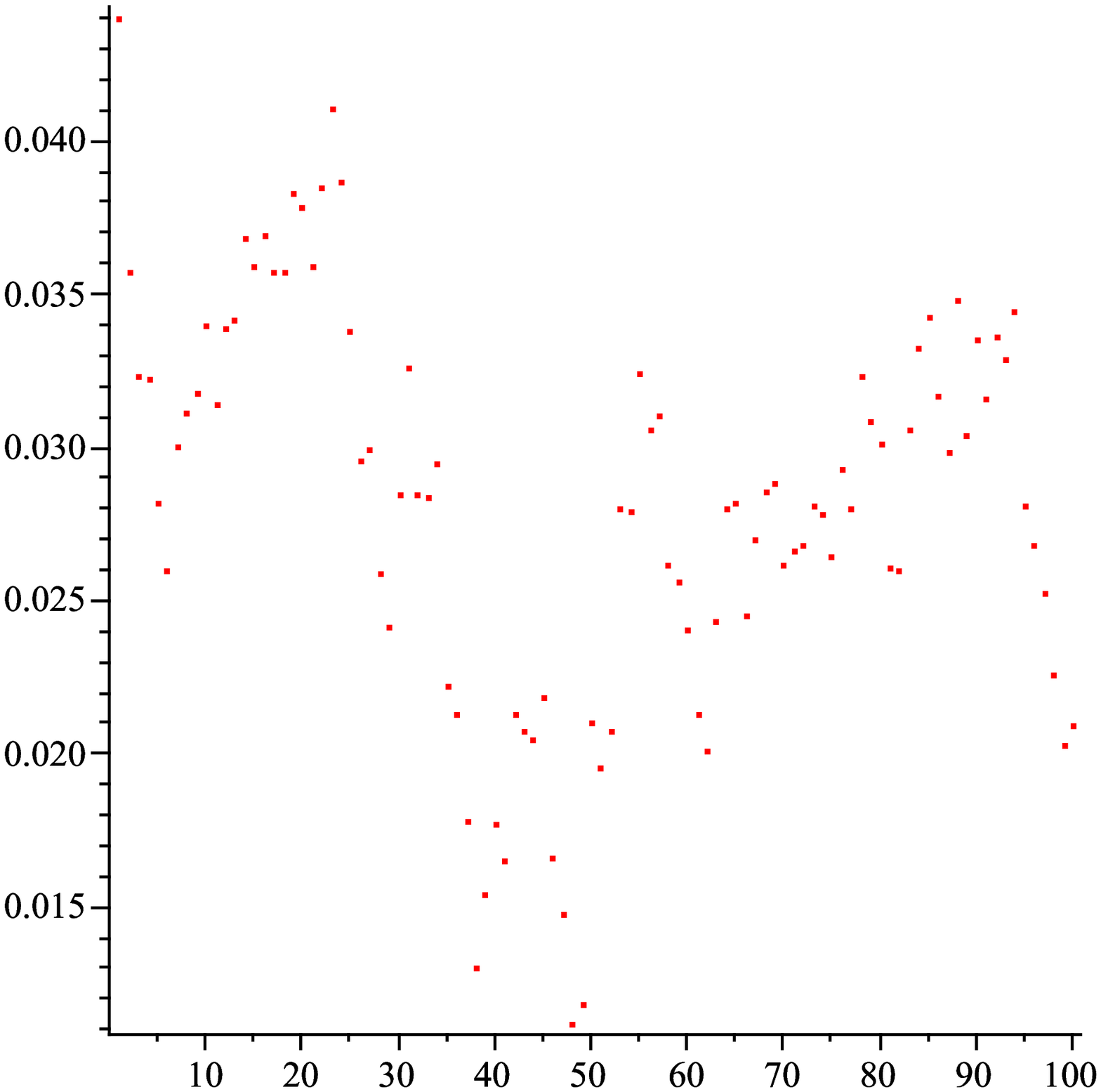}& 
\epsfxsize=5cm \epsfysize=4cm
\epsffile{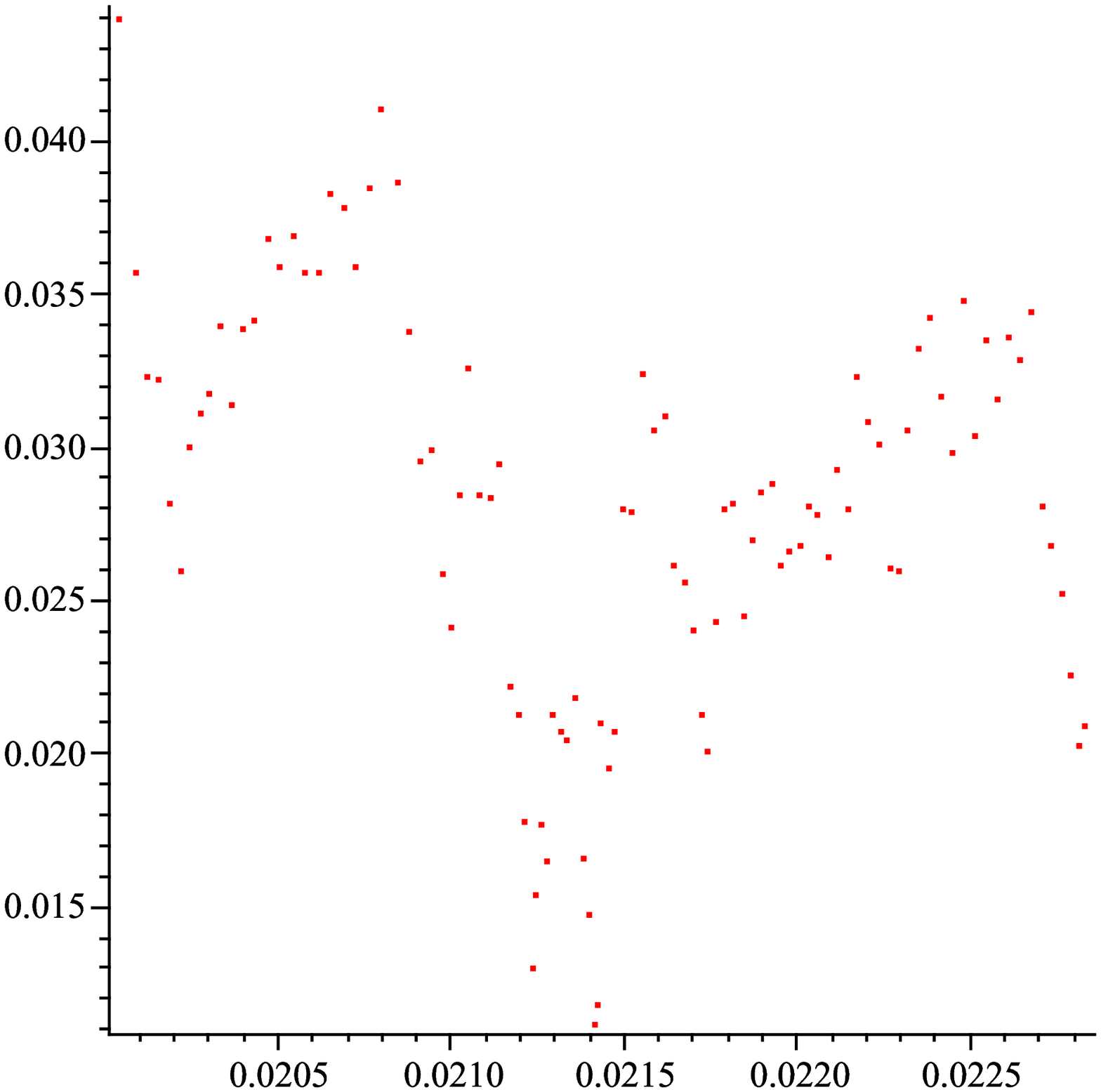}\\
\textit{ Fig16:} $(n,q(n,\omega))$& \textit {Fig17:} $(n,p(n,\omega))$&
\textit{ Fig18:} $(q(n,\omega),p(n,\omega))$\\
    \end{tabular}
\end{center}

\section{CONCLUSIONS}

In this paper we present generalization of fractional Riemann-Liouville integral, Wiener process,  and we have defined the generalized fractional stochastic, Liu and hybrid equations.  The mixture between generalized fractional Wiener process and generalized fractional Liu process results as the generalization of fractional hybrid differential equations.  We defined generalized fractional hybrid HP equations and generalized fractional hybrid Hamiltonian equations. The first order Euler scheme is presented and implemented for particular parameters. 
In the future work, we will consider other problems that deal with
stochastic fractional HP principle.

\end{document}